\documentclass[9pt,journal,twocolumn]{IEEEtran}
\usepackage[noadjust]{cite}
\usepackage{longtable}
\usepackage{amsmath}
\usepackage{graphicx, amssymb}
\usepackage{epstopdf}
\usepackage{balance}
\usepackage[dvips]{epsfig}
\usepackage{color}
\usepackage{comment}
\usepackage{setspace}
\usepackage[dvipsnames]{xcolor}
\colorlet{LightRubineRed}{RubineRed!70!}
\colorlet{Mycolor1}{green!10!orange!90!}
\definecolor{Mycolor2}{HTML}{00F9DE}

\newtheorem{thm}{Theorem}

\newtheorem{rem}[thm]{Remark}

\newtheorem{lem}[thm]{Lemma}
\newtheorem{cor}[thm]{Corollary}

\newtheorem{pro}[thm]{Proposition}

\def\L2{{\cal L}_2}

\newcommand{\beq}{\begin{equation}}
\newcommand{\eeq}{\end{equation}}

\newcommand{\carre} {\hfill $\blacksquare$}

\def\BibTeX{{\rm B\kern-.05em{\sc i\kern-.025em b}\kern-.08em
    T\kern-.1667em\lower.7ex\hbox{E}\kern-.125emX}}

\begin{document}
%\begin{frontmatter}

\title{Laplacian Dynamics on Cographs:\\
Controllability Analysis through Joins and Unions} %!PN
%Towards Modular Controllability Analysis on Networks} %!PN

\author{Shima Sadat Mousavi, Mohammad Haeri, and Mehran Mesbahi% <-this % stops a space
\thanks{S.S. Mousavi was with the Department of Electrical Engineering, Sharif University of Technology, Tehran, Iran. She is now with the Department of Civil, Environmental and Geomatic Engineering, ETH Z\"urich, 8093 Z\"urich, Switzerland (e-mail: shimasadat$\_$mousavi@ee.sharif.edu, mousavis@ethz.ch).}
\thanks{M. Haeri is with the Department of Electrical Engineering, Sharif University of Technology, Tehran, Iran (e-mail: haeri@sharif.ir).}%
\thanks{M. Mesbahi is with the Department of Aeronautics and Astronautics, University of Washington, WA 98195 (e-mail: mesbahi@uw.edu).}
%\thanks{The paper was presented in part at the 55th IEEE Conference on Decision and Control, Las Vegas, December 10-13, 2016 (see \cite{mousavi2016controllability}).}
}

\maketitle
%\title{Using the Style File IEEEtran.sty} 

%\thispagestyle{plain}\pagestyle{plain}

\begin{abstract}
In this paper, we examine the controllability of Laplacian dynamic networks on cographs. % dynamics follow a consensus-type coordination protocol, or the so-called Laplacian dynamics, and their 
Cographs appear in modeling a wide range of networks and include as special instances, the threshold graphs. In this work, we present  necessary and sufficient conditions for the controllability of cographs, and provide an efficient method for selecting a minimal set of input nodes from which the network is controllable. In particular, we define a sibling partition in a cograph and show that the network is controllable if all nodes of any cell of this partition except one are chosen as control nodes. The key ingredient for such characterizations is the intricate connection between the modularity of cographs and their modal properties.
Finally, we use these results to characterize the controllability conditions for certain subclasses of cographs.       
\end{abstract}


%\end{frontmatter}

%-------------------------------------------------------------------------------

\section{Introduction}

Networks have become the backbone of the modern society.
%Social networks, the Internet, and energy systems are examples of critical networks that we rely on their operation in our daily lives.
%As such, the control, security, and management of %these and other types of 
%networks are of paramount importance, providing a rich class of system theoretic questions for the control community~\cite{mesbahi2010graph}.
One foundational class of questions on networked systems
pertains to their controllability~\cite{tanner2004controllability,mousavi2018structural,mousavi2017robust,
mousavi2018null}. While there are classical tests to check the controllability of linear time-invariant (LTI) systems, their application to large-scale networks is numerically infeasible. Moreover, finding a minimum cardinality set of input nodes  ensuring the controllability of a  network is NP-hard. To overcome these issues, an alternative set of approaches involves adopting graph-theoretic techniques and connecting the controllability of a network to its topological features.
  %
 % Such an approach to the network controllability problem can also provide a framework for designing topologies that have favorable system theoretic properties. 
  In this direction,  controllability analysis of networks with the so-called Laplacian dynamics has gained a lot of attention, partially due to their relevance in distributed algorithms  as consensus, distributed estimation, and nonlinear synchronization~\cite{tanner2004controllability, zhang2014upper, ji2015protocols,yaziciouglu2016graph,ji2017new}.
  \color{black}
  In its most basic form, this dynamics is realized when the states of a network follow  a consensus-type coordination protocol \cite{olfati2007consensus}. In such a dynamics, a subset of nodes in the network-known as leaders-are assumed to be controlled by external commands, while the other nodes-referred to as followers-follow the consensus (nearest-neighbor interaction) protocol. The controllability analysis of  leader-follower Laplacian networks is of great interest in scenarios such as formation control, human-swarm interaction, and network security~\cite{ji2009interconnection,mesbahi2010graph}. 
  \color{black}
    %exploring network controllability from the topology perspective.  %investigating the controllability of a large scale network and finding a set of control nodes with the minimum size which  through the classical rank conditions is computationally impossible. An alternative way is exploring network controllability from its topology perspective. %Many works in this field are focused on networks following a Laplacian dynamics. In fact, these  %of a network through the features of its topology which is described by a graph. %  finding  from technological ones to real-world systems, constitute an important part of Controllability analysis of linear time-invariant (LTI) networks through features of their corresponding graph has been an active area of research in the system and control community for the past several decades. 
%

There are two classes of results in the literature on the controllability analysis of Laplacian networks. In the first setting, 
%regardless of the topology of  network, 
necessary or sufficient conditions have been provided for network controllability. These conditions have been mainly stated in terms of notions such as graph symmetry \cite{rahmani2009controllability,chapman2015state}, equitable partitions \cite{rahmani2009controllability, egerstedt2012interacting,zhang2011controllability, zhang2014upper, cao2013class,aguilar2017almost}, distance partitions \cite{zhang2011controllability, zhang2014upper}, and pseudo monotonically increasing sequences \cite{yazicioglu2012tight, yaziciouglu2016graph}. For example, the existence of a symmetry with respect to the leaders or control nodes of a network is known to be a sufficient condition for its uncontrollability \cite{rahmani2009controllability}. There are, however, drawbacks to this line of work for analyzing large-scale networks. First, the known graph-theoretic conditions are not necessary {\em and} sufficient for network controllability; rather, these conditions are often used to obtain lower or/and upper bounds on the dimension of the controllable subspace. Furthermore, most of these results cannot be utilized for efficiently selecting input nodes ensuring the controllability of the network. For instance,  finding a minimum cardinality set of nodes breaking symmetries for general networks is NP-hard \cite{chapman2015state}. We also mention the results reported in~\cite{ji2017new}, where by identifying the structure of controllability destructive nodes, necessary and sufficient controllability conditions for Laplacian networks of size five or less have been established. \color{black}

%The results in the literature on the controllability analysis of networks with Laplacian dynamics can be classified into two categories. In the first category,  a general topology has been  considered  for the network, and  some necessary \emph{or} sufficient conditions for its controllability from a graph-theoretic point of view have been  presented. These conditions have been stated in terms of notions such as graph symmetry \cite{rahmani2009controllability,chapman2015state}, equitable partitions \cite{rahmani2009controllability, egerstedt2012interacting,martini2010controllability,zhang2011controllability, zhang2014upper, cao2013class,aguilar2017almost}, distance partitions \cite{zhang2011controllability, zhang2014upper}, and pseudo monotonically increasing sequences \cite{yazicioglu2012tight, yaziciouglu2016graph}. However, they have a few important shortcomings. In fact, none of these conditions are necessary \emph{and} sufficient for the network controllability; rather, they result in deriving only some lower or/and upper bounds on the dimension of the controllable subspace. More importantly, these results cannot be utilized for efficient selection of control nodes rendering a network controllable. For example, it is known that the existence of a symmetry in the structure of a network with respect to its control nodes is destructive to its controllability \cite{rahmani2009controllability}, but finding a minimum cardinality set of nodes breaking all symmetries in a general network is NP-hard \cite{chapman2015state}.  

In order to derive stronger and readily applicable network-centric controllability conditions, in the second class of results, Laplacian networks with special graph topologies have been considered.  %The second class of results includes the works considering a special graph topology for the Laplacian network in order to derive stronger and more accurate controllability conditions for this graph.
 In this case, controllability of networks with embedded path graphs \cite{parlangeli2012reachability, mousavi2016controllability}, cycle graphs \cite{parlangeli2012reachability},  \cite{liu2018controllability}, \color{black} complete graphs \cite{zhang2011controllability}, circulant graphs \cite{nabi2013controllability}, multi-chain graphs \cite{hsu2017necessary}, grid graphs \cite{notarstefano2013controllability}, and tree graphs \cite{ji2012leaders} have been investigated. These approaches rely on the pattern of the Laplacian eigenvectors in conjunction with the Popov-Belevitch-Hautus (PBH) test to facilitate the controllability analysis. Moreover, a complete characterization of the eigenspaces of these graphs leads to efficient procedures for selecting the minimum number of control nodes from which the network is controllable. By considering the different methods of combining or growing controllable networks,  this class of results can  be also applied when designing network structures with desired controllability properties \cite{ tran2018generalized, ji2016design}. \color{black}However, the class of Laplacian networks with efficient graph-theoretic controllability conditions is still limited. In this paper, we further expand the applicability of such graph-theoretic conditions by examining the controllability of Laplacian networks defined over \emph{cographs}.

%The second category includes the works that consider a  special class of graphs and  study the controllability of networks defined on these graphs \cite{aguilar2015graph}. For example, controllability of networks with path graphs \cite{parlangeli2012reachability, mousavi2016controllability}, cycle graphs \cite{parlangeli2012reachability}, complete graphs \cite{zhang2011controllability}, circulant graphs \cite{nabi2013controllability}, multi-chain graphs \cite{hsu2017necessary}, grid graphs \cite{notarstefano2013controllability}, and tree graphs \cite{ji2012leaders} have already been explored, and some stronger conditions for their controllability have been derived. In particular, for some of these graphs, the minimum number of control nodes from which the associated network is controllable has been determined. Note that the stronger controllability conditions derived for these special classes of graphs are resulted from a better characterization of the eigenvectors associated with their Laplacian. In fact, based on the Popov-Belevitch-Hautus (PBH) test, the controllability of a system depends completely on its associated eigenvectors. Then, by identifying the eigenspace of a network (i.e., the space of eigenvectors associated with each eigenvalue of the Laplacian matrix), the controllability problem can be totally solved.

%Adopting a similar approach, in this paper, we consider the controllability problem for the Laplacian networks defined on cographs. 
Cographs have been independently rediscovered and reintroduced by different authors; as such, they assume multiple equivalent definitions. 
For example, in such graphs, there is no induced subgraph isomorphic to a path of size four, and accordingly they are called $P_4$-free graphs~\cite{corneil1981complement}. Moreover, some authors refer to cographs as decomposable graphs \cite{merris1998laplacian}, or complement-reducible graphs \cite{corneil1981complement}, due to the fact that  they can be generated through recursive operations of \emph{joins} and \emph{unions} starting from isolated nodes~\cite{biyikoglu2007laplacian}.  The sequence of these operations  leads to a unique rooted tree representation of a cograph,  referred to as a cotree \cite{corneil1981complement}.  Cographs arise in disperate areas of computer science and mathematics and find applications in areas such as scheduling~\cite{corneil1984cographs,corneil1985linear} and orthology detection~\cite{hellmuth2013orthology}\color{black}. In fact, thanks to their structural properties, many algorithmic problems that are NP-hard for general networks can be solved in a polynomial time over cographs~\cite{habib2005simple}. Cographs have  a close relationship with series-parallel networks that are used to model biological and electrical systems \cite{cheng2010efficient,corneil1981complement}. Furthermore, there has been interest in identifying cograph communities and functional modules in social and biological networks in order to better reveal their local and global structures and functions~\cite{jia2015defining,jia2018viewing}. \color{black} Cographs include other known classes of graphs, including complete graphs, complete bipartite graphs, cluster graphs, Turan graphs, and trivially-perfect graphs.  In particular, threshold graphs are an important subclass of cographs with numerous applications in modeling social and psychological networks and synchronizing parallel processes \cite{saha2014intergroup, mahadev1995threshold}. \color{black} %There are some different representations for a threshold graph as well; for instance, it can be uniquely determined by a binary construction sequence \cite{hagberg2006designing}.  
 In the meantime, a connected threshold graph has certain limitations in modeling networked control systems; for example, it has at least one dominating node (that is adjacent to all other nodes), while a general cograph might not have such a node. 
 
\color{black}

In \cite{aguilar2015laplacian}, the controllability of threshold graphs from a single control node has been explored. % Classifying the single control input networks into three groups, namely essentially controllable, conditionally controllable, and completely uncontrollable graphs, the work~\cite{aguilar2015laplacian} characterizes necessary and sufficient conditions for a threshold graph to be completely uncontrollable. %It has also been proven that a threshold graph of size $n$  is controllable from a single control input if and only if it is an anti-regular graph with $n-1$ different degrees. 
Subsequently, in \cite{hsu2016controllability}, the results of \cite{aguilar2015laplacian} have been extended to multi-input networks; however, the results provided in  \cite{hsu2016controllability} are restrictive in the sense that it examines threshold graphs with only one repeated degree. In \cite{mousavi2018controllability}, the controllability problem of a general threshold graph has been solved for the case in which any input signal is assumed to be injected into only one node.  Furthermore, in \cite{hsu2019minimal}, the results of \cite{mousavi2018controllability} have been  extended to the case where the input matrix has binary entries. However, we note that unlike a general cograph, in a threshold graph, information about the set of eigenvalues and eigenvectors  can be inferred from the sequence of node degrees. As such, the results of \cite{mousavi2018controllability} and \cite{hsu2019minimal} cannot readily be applied in a more general setting, namely for networks that are characterized by cographs. In fact, the analysis approach adopted in the aforementioned works--relying on the creation sequence and node degrees of a threshold graph--does not apply for a general cograph.    \color{black}  In this paper, we take a step towards studying the controllability problem for Laplacian networks defined on  cographs. 
%To the best of our knowledge, this is the first study on the controllability of a general cograph. %The main contributions of the paper are as follows: First,

 Our first contribution is providing the spectrum and an associated modal matrix for a cograph. This is accomplished by considering the cotree representation of cographs, and subsequently showing that the set of nontrivial eigenvalues (respectively, eigenvectors) of a cograph is an updated version of the nontrivial eigenvalues (respectively, eigenvectors) generated at each internal node of the associated cotree. In this direction, we also illustrate some properties of eigenvalues of a cograph and their eigenspaces, using the structural feature of the associated cotree; this is discussed in \S\ref{cograph-spectra}.  %we first characterize the Laplacian eigenspace of a cograph.  
 \color{black}
 As the second and main contribution of this work, we establish necessary and sufficient conditions for controllability of Laplacian networks on cographs in \S\ref{controllability}. In this direction, based on the fact that the uncontrollability of a network results from the zero entries of its eigenvectors, we identify all (and the only) nodes rendering a cograph controllable.  In fact, we decompose a cograph into structurally equivalent subgroups or cells, essentially playing similar roles in the network dynamics.
By defining a sibling partition in a cograph, we then demonstrate that these cells include sibling nodes that interact similarly with all other nodes in the graph. %are structurally similar to the rest of nodes. 
Thus, in order  to break ``structural symmetries" in a cograph, all nodes of any cell except one should be directly controlled. Particularly, it is proven that the minimum number of control nodes to completely control a cograph is the difference between its size and the number of cells of its sibling partition. Finally, we provide an alternate  approach for obtaining the results reported in~\cite{mousavi2018controllability,hsu2019minimal}, when the controllability conditions for general cographs are interpreted in the context of its subclasses, such as the threshold graphs.

 \section{Preliminaries} 
  
%\begin{deff}  The set of special  orthogonal matrices is defined by $\mbox{SO}(n) = \{ X\in \Rdd n n: X^\top X =I, \det(X) = 1 \}$.\end{deff}
In this section, necessary preliminaries for our
subsequent discussion are reviewed.
\newline
\emph{Notation:} The set of real numbers is denoted by $\mathbb{R}$. For a set $\mathcal{S}$, its cardinality is denoted by $|\mathcal{S}|$. For a matrix $M\in\mathbb{R}^{p\times q}$ and a set of indices $\mathcal{S}$, $M_{\mathcal{S},:}\in\mathcal{R}^{|\mathcal{S}|\times q}$ is a submatrix of $M$ whose rows are the indices from $\mathcal{S}$. %$M_{ij}$ is the entry of $M$ in its $i$th row and $j$th column, and $M_{i,:}$ and $M_{:,j}$ represent the $i$th row and $j$th column of $M$. For two sets of indices $s_1$ and $s_2$, $M_{s_1,s_2}$ is a submatrix of $M$ whose rows (respectively, columns) are the indices from $s_1$ (respectively, $s_2$). 
The $n\times n$ identity matrix is denoted by $I_n$, and $e_j$ represents its $j$th column. The vectors of all 1's  and all 0's with size $n$ are respectively denoted by  $\textbf{1}_n$ and $\textbf{0}_n$. The $n\times m$ matrix of all 1's (respectively, all 0's) is designated as $\textbf{1}_{n\times m}$ (respectively, $\textbf{0}_{n\times m}$). For notational convenience, for a vector $v\in \mathbb{R}^n$ and a scalar $m\in \mathbb{R}$, we write $v+m$ to represent $v+m\textbf{1}_{n}$.

%\subsection{Preliminaries on Graphs and Patterned Matrices}

\textit{Graph:} A directed graph\footnote{All graphs in this paper are assumed to be unweighted, simple, and loop-free.} 
$G$ of size $n$ is represented by $G=(V,E)$, where 
$V=\{1,\ldots,n\}$ is its node set, and $E\subset V\times V$ denotes its edge set. % defined by $E(G)=\{\{i,j\}: i,j\in V(G), \:\:i\neq j,:\: \mbox{there}\: \:\:\mbox{is}\: \:\:\mbox{an}\: \:\: \mbox{edge}\: \:\: \mbox{between}\:\:\: i\: \: \: \mbox{and} \:\:\: j\}$.
We say $(i,j)\in E$ if there is a directed edge from the node $i$ to the node $j$. A directed path from the node $i_1\in V$ to the node $i_k\in V$ is a sequence of distinct nodes $(i_1,i_2,\ldots,i_k)$, where for every $1\leq j\leq k-1$, $(i_j,i_{j+1})\in E$. The graph $G$ is undirected if for every edge $(i,j)\in E$, we have $(j,i)\in E$; in this case, we write $\{i,j\}\in E$, and we refer to node $j$ (respectively, $i$) as the neighbor of the node $i$ (respectively, $j$). For an undirected graph $G$, we denote by $N(i)$ the set of neighbors of $i\in V$. The degree of the node $i$ is defined as $d(i)=|N(i)|$. The degree matrix of an undirected graph $G$ is defined as $\Delta(G)=\mathrm{diag}(d(1),\ldots,d(n))$. The corresponding Laplacian matrix $L(G)$ is given by $L(G)=\Delta(G)-A(G)$, where $A(G)$ is the (0,1)-adjacency matrix associated with  $G$. %Let the nodes of the graph be ordered and indexed in a way that $d_i\leq d_{i+1}$, for $1\leq i\leq n-1$. 
A \emph{complete} graph $G=(V,E)$ is an undirected graph such that for all $i,j\in V$, $i\neq j$, $\{i,j\}\in E$; it is denoted by $K_n$. Consider two disjoint sets $V_1$ and $V_2$ of respectively size $n_1$ and $n_2$ such that $V=V_1\cup V_2$. A \emph{complete bipartite} graph $G=(V,E)$, denoted by $K_{n_1,n_2}$, is an undirected graph such that for any pair of nodes $i,j\in V_k$, $k=1,2$, $\{i,j\}\notin E$, while for any $i\in V_1$ and $j\in V_2$, $\{i,j\}\in E$. A \emph{path} graph of size $n$, denoted by $P_n$, is a graph whose nodes can be indexed by $1,2,\ldots, n$ in such a way that for all $1\leq i\leq n-1$, $\{i,j\}\in E$ if and only if $j=i+1$. If two nodes of degree one in $P_n$ are connected, a \emph{cycle} graph $C_n$ is obtained.  

\emph{Rooted trees:} %An undirected tree is generated from a single node by repeatedly adding a node which is adjacent to only one of the old nodes. Now, assign a direction to every edge of a tree graph. If there is a special node $r\in V$, called the root, and there is a unique directed path from every node of the tree to $r$, then it is a rooted tree. % Now, consider a sequence of distinct nodes $(i_1,i_2,\ldots,i_k)$ associated with a directed path from $i_1$ to $i_k$. We define $P_{i_1}^{i_k}=\{(i_1,i_2),\ldots, (i_{k-1},i_k)\}$ as a collection of edges defining a path from $i_1$ to $i_k$. 
Consider an undirected tree graph, and assign a direction to any of its edges. The new directed graph is a rooted tree and denoted by $\mathcal{T}=(V^{\mathcal{T}},E^{\mathcal{T}})$ if for a special node $r\in V^{\mathcal{T}}$, called the root, there is a unique directed path from $r$ to every node of $\mathcal{T}$.   \color{black}
%A rooted tree denoted by $\mathcal{T}=(V^{\mathcal{T}},E^{\mathcal{T}})$ is a directed graph which has a node $r\in V^{\mathcal{T}}$, called the root, so that there is a unique directed path from every node of $\mathcal{T}$ to $r$. 
% Now, consider a sequence of distinct nodes $(i_1,i_2,\ldots,i_k)$ associated with a directed path from $i_1$ to $i_k$. We define $P_{i_1}^{i_k}=\{(i_1,i_2),\ldots, (i_{k-1},i_k)\}$ as a collection of edges defining a path from $i_1$ to $i_k$. 
For a node $i\in V^{\mathcal{T}} $, a node $j\in V^{\mathcal{T}}$ such that $(i,j)\in E^{\mathcal{T}}$ (respectively, $(j,i)\in E^{\mathcal{T}}$) is called a \emph{child} (respectively, \emph{parent}) of $i$. %A group of nodes with the same parent is referred to as \emph{siblings}.
 A node $j\in V^{\mathcal{T}}$ is called a \emph{descendant} (respectively, \emph{ancestor}) of node $i\in V^{\mathcal{T}}$ if there is a directed path from $i$ to $j$ (respectively, from $j$ to $i$).  The \emph{lowest common ancestor} of two nodes $k,l\in V^{\mathcal{T}}$ is the shared ancestor of $k$ and $l$, which is located farthest from the root $r$ of $\mathcal{T}$. \color{black} A node $i$ is called a \emph{leaf} if it has no child; otherwise, it is an \emph{internal node} of $\mathcal{T}$. The set of children of an internal node $v$ is given by ${\mathcal{C}}(v)$, and its size is denoted by $c(v)$. Moreover, the set of leaves descending from the internal node $v$ is represented by $\mathcal{L}(v)$; we define $l(v)=|\mathcal{L}(v)|$. The unique path from a node $w$ to its descendant $v$ is given by $\mathcal{P}_v^w$.  A group of  leaves \color{black}  with the same parent in the rooted tree $\mathcal{T}$ is referred to as \emph{siblings}. %leaf $v$ to the root of  $\mathcal{T}$ is denoted by $P_{\mathcal{T}}(v)$.

\emph{Example:} In the rooted tree depicted in Fig. \ref{cograph} (b), any node $i$, $1\leq i\leq 8$, is a leaf; while each node $v^j$, $1\leq j\leq 5$, is an internal node.    Node $v^3$ is the root, and we have $\mathcal{C}(v^3)=\{4,v^2, v^4\}$. Thus, $c(v^3)=3$. Also, one has $\mathcal{L}(v^4)=\{5,6,7,8\}$, and $l(v^4)=4$. Leaves 6, 7, and 8 have the same parent $v^5$, and thus are siblings. The path $\mathcal{P}_{v^1}^{v^3}$ can be described by the sequence $(v^3,v^2,v^1)$.  The lowest common ancestor of leaves 5 and 8 is $v^4$. \color{black}
  
  \color{black}
 \emph{Eigenpairs:} Consider an undirected graph $G=(V,E)$. For notational convenience, by eigenvalues and eigenvectors of $G$, we mean the eigenvalues and eigenvectors of its Laplacian matrix $L(G)$. Since $L(G)$ (for an undirected graph $G$) is symmetric and nonnegative, all of its eigenvalues are real and nonnegative \cite{godsil2013algebraic}\color{black}. Moreover, its smallest eigenvalue is zero with the associated eigenvector $\textbf{1}_n$.  
The vector $(0,\textbf{1}_n)$ is known as a trivial eigenpair for any undirected graph $G$.
  Now, let $\Lambda(G)=(\lambda_2,\ldots,\lambda_n)$ be the nontrivial spectrum of $G$ including its nontrivial eigenvalues,  %Define $\bar{\Lambda}=(\lambda_2,\ldots, \lambda_n)$ as the nontrivial spectrum of $G$, 
and note that if $G$ is connected, $0\notin {\Lambda}(G)$. Next, let $\nu_i\in \mathbb{R}^n$ be a nonzero eigenvector of $G$ associated with $\lambda_i$, where $L(G)\nu_i=\lambda_i\nu_i$. Then, we define ${\mathcal{V}}(G)=[ \nu_2,\ldots,\nu_n]$ as a full rank  \color{black} nontrivial modal matrix of $G$ associated with ${\Lambda}(G)$. %, where $\bar{\mathcal{V}}=[\nu_2,\ldots,\nu_n]\in\mathbb{R}^{n\times (n-1)}$. % Since $L(G)$ is symmetric and $L(G)\geq0$ (for an undirected graph $G$), all of its eigenvalues are real and nonnegative. It is known that $\mbox{min}(\Lambda(G))=\mbox{min}(\lambda_1,\ldots,\lambda_n)=0$ and its associated eigenvector is $\textbf{1}_n$. If $G$ is a connected graph, $0$ is a simple eigenvalue of $G$. % Let $\Lambda(G)=(\lambda_1,\ldots,\lambda_n)$ be the spectrum of the graph $G$, where $ \lambda_1\leq \lambda_2\leq\ldots\leq \lambda_n$. Then, $\lambda_1=0$, and if $G$ is connected, we have $\lambda_2\neq 0$.
Let $\tilde{\lambda}_1$, $\ldots$, $\tilde{\lambda}_r$ be the $r$ distinct eigenvalues in ${\Lambda}(G)$.  We can rewrite the nontrivial spectrum of a connected $G$ as  ${\Lambda}(G)=(\tilde{\lambda}_1^{(q_1)},\ldots,\tilde{\lambda}_r^{(q_r)})$,  where $q_i$ is the algebraic multiplicity of the nonzero eigenvalue $\tilde{\lambda}_i$.  %Then, $\mathcal{M}=\mbox{max}\{q_1,\ldots, q_r\}$ is the maximum multiplicity of eigenvalues of $G$. 
 Since $L(G)$ is symmetric, for an eigenvalue $\tilde{\lambda}_i$ with the multiplicity $q_i$, there are $q_i$ independent eigenvectors  spanning the eigenspace associated with $\tilde{\lambda}_i$. Let $\mathcal{V}^{(i)}\in \mathbb{R}^{n\times q_i}$ be a full rank matrix where $L(G)\mathcal{V}^{(i)}=\tilde{\lambda}_i\mathcal{V}^{(i)}$, $i=1,\ldots,r$. Then, the nontrivial modal matrix associated with ${\Lambda}(G)$ for a connected $G$ can be written as ${\mathcal{V}}(G)=[\mathcal{V}^{(1)},\ldots,\mathcal{V}^{(r)}]$.  

 \subsection{Cographs}

%\emph{Cographs:} 
In this part,
 %the notion of
 cographs and related concepts are reviewed.

Let $G_1=(V_1 ,E_1)$ and $G_2=(V_2,E_2)$ be two disjoint undirected graphs of respectively, size $n_1$ and $n_2$.  A graph $G=(V,E)$ is the \emph{union} of  $G_1$ and $G_2$ if $V=V_1\cup V_2$, and $E=E_1\cup E_2$; such a graph $G$ is written as $G=G_1+G_2$. A graph $G=(V,E)$ is the \emph{join} of $G_1$ and $G_2$ if $V=V_1\cup V_2$, and $E=E_1\cup E_2 \cup \{\{i,j\}: i\in V_1 \:\mbox{and}\: j\in V_2 \}$; thus $G$ is represented by $G=G_1*G_2$.  Join and union operations obey the commutative and transitive properties. For example, we have $G_1+G_2=G_2+G_1$, and $G_1*(G_2*G_3)=(G_1*G_2)*G_3$. However, they are not distributive.
\color{black}
%The \emph{union} of the two graphs is a graph of size $n=n_1+n_2$, which is defined as $G_1+G_2=(V_1\cup V_2, E_1\cup E_2)$. Moreover, the \emph{join} of the two graphs represented by $G_1 * G_2$ is obtained from $G_1+G_2$ by adding new edges from each node of $G_1$ to any node of $G_2$. 

A graph is called a \emph{cograph}\footnote{In this paper, cographs are assumed to be undirected.} (or a \emph{decomposable graph}) if it can be constructed from isolated nodes by recursively performing join and union operations. More formally,  a graph with a single node (i.e., $K_1$) is a cograph, and if $G_1,\ldots,G_k$, for some $k>1$, are cographs, then $G_1+\ldots +G_k$ and $G_1*\ldots* G_k$ are cographs as well. Also, based on an equivalent definition, a graph is a cograph if it has no induced subgraphs isomorphic to $P_4$ \cite{hammer1996laplacian}. Thus, for example, the cycle graph $C_5$ is not a cograph, while the complete graph $K_4$ is. 
\color{black}

A cotree $\mathcal{T}=(V^{\mathcal{T}},E^{\mathcal{T}})$ associated with a connected cograph $G=(V,E)$ is a rooted tree  whose leaves correspond to the nodes of the cograph. We index the leaves of a cotree by $i$, where $1\leq i\leq n$; while the internal nodes are represented here by $v^j$, for $1\leq j\leq h$. \color{black}  \color{black}   The root of the cotree $r$ is labeled as 1, and its internal nodes  are labeled  0 or 1.  For an internal node $v^k$, $\mathrm{lab}(v^k)$ provides the label of $v^k$. Let  $\mathcal{T}_{(z)}$ be a subtree of $\mathcal{T}$ which is rooted at some node $z\in V^{\mathcal{T}}$. % For any node $z\in V^{\mathcal{T}}$, 
Then, $\mathcal{T}_{(z)}$ corresponds to an induced subgraph of $G$ defined on the leaves which are descendants of $z$. We denote this subgraph by $G_{(z)}$, and call it a cograph associated with $z$. % $G_{(z)}$ is itself a cograph.
 If $z$ is a leaf of $\mathcal{T}$, $G_{(z)}=(\{z\},\emptyset)$. % is  a graph of the single node $z$. 
In addition, if $z$ is an internal node that is labeled as 0 (respectively, 1), $G_{(z)}$ is the union (respectively, join) of the cographs associated with the children of $z$ \cite{corneil1981complement}.  
 %

  %On the other hand, if $z$ is labeled 1,  $G_{(z)}$ is a join of subgraphs corresponding to  the children of $z$ . 
%$\mathcal{T}_z$ is associated with a graph of the single node $z$ in $G$. In addition, if $z$ is an internal node of $\mathcal{T}$ that is labeled 0, $\mathcal{T}_z$ corresponds to a subgraph of $G$ which is the union of subgraphs associated with the children of $z$. On the other hand, if $z$ is labeled 1, the corresponding subgraph of $\mathcal{T}_z$ in $G$ is a join of subgraphs corresponding to  the children of $z$ \cite{corneil1981complement}. 

 %Any connected cograph $G$ can be represented by a  cotree $\mathcal{T}$ whose root is labeled 1.
 Any cograph  $G$ can be represented by a cotree $\mathcal{T}$, and if for any leaf $v$ of $\mathcal{T}$, the labels on the internal nodes of the path $\mathcal{P}_v^r$ alternate between 0 and 1, this representation is unique. A cograph $G=(V,E)$ can be recognized in $\mathcal{O}(|V|+|E|)$, while its associated cotree can be built in the same time-complexity~\cite{corneil1985linear}. \color{black} %Note that a cograph $G$ is connected if and only if the root of $\mathcal{T}$ is labeled 1. 
 Alternatively, one can form a cograph $G$ from a given cotree $\mathcal{T}$. In this direction, two nodes $i,j\in V$ are neighbors in $G$ if and only if the lowest common ancestor of the leaves $i,j\in V^{\mathcal{T}}$ is labeled 1. %This node is called the \emph{lowest common ancestor} of  leaves $i$ and $j$. 

In a cograph $G=(V,E)$, two nodes $i,j\in V$ are called siblings if the leaves $i$ and $j$ in the corresponding cotree are siblings. By this definition, it is known that $i,j\in V$ are siblings if $N(i)\setminus \{j\}=N(j)\setminus \{i\}$ \cite{corneil1981complement}. 

\begin{figure}[hbt]
\includegraphics[width=.4\textwidth]{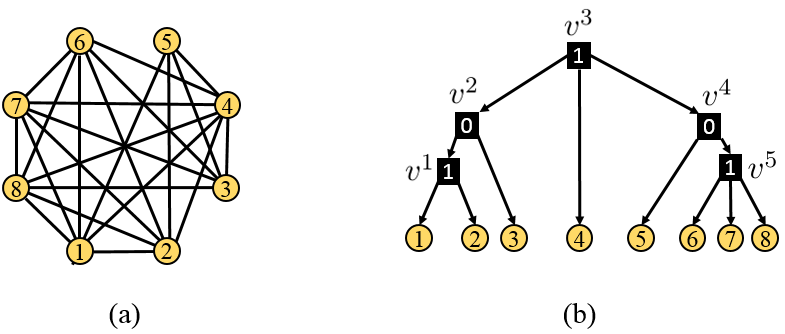}
\centering
\caption{a) Cograph $G$, b) Associated cotree $\mathcal{T}$.}
\label{cograph}
\end{figure}

\emph{Example:} In Fig. \ref{cograph},  a cograph along with its associated cotree are illustrated. The cograph $G$ can be constructed through successive joins and unions as $G=[(K_1*K_1)+K_1]*K_1*[K_1+(K_1*K_1*K_1)]$. One can see that nodes 1, 2 and nodes 6, 7, 8 are siblings in this cograph. We have $\mathrm{lab}(v^1)=1$, and $\mathrm{lab}(v^4)=0$. The graph $G_{(v^2)}$ is an induced subgraph of $G$ defined on nodes $1,2,3$. Since $\mathrm{lab}(v^2)=0$, $G_{(v^2)}$ is a union of $G_{(v^1)}$ and $G_{(3)}$.  The lowest common ancestor of nodes 5 and 8 is $v^4$, that is labeled 0; thus, these two nodes are not neighbors in $G$.
 \color{black}

\subsection{Threshold Graphs}

By starting from a single node, a threshold graph that is a special subclass of cographs is constructed by repeatedly adding a single node to the old graph through the join or the union operation. In other words, $K_1$ is a threshold graph; and if $G'$ is a threshold graph, $G'+K_1$ and $G'*K_1$ are threshold graphs as well.  % Let us denote by $v^i$ the $i$th farthest internal node from the root of the cotree associated with a threshold graph. Also, index as $j$ the $j$th node added to the old graph through join or union. 
%
 %One can associate a unique binary \emph{construction sequence} $T^G\in\{0,1\}^n$ to a threshold graph $G$ of size $n$, where $T^G(1)=0$, and for $1<i\leq n$, $T^G(i)=0$ (respectively, $T^G(i)=1$) if the $i$th node is added to the former graph through the union (respectively, join) operation \cite{hagberg2006designing}. %Note that $T^G$%In fact, any threshold graph $G$ can be uniquely determined by its construction sequence $T^G$. 
 
\begin{figure}[hbt]
\includegraphics[width=.4\textwidth]{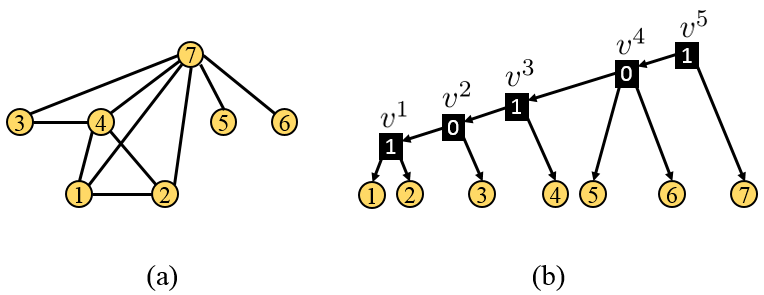}
\centering
\caption{a) Threshold graph $G$, b) Associated cotree $\mathcal{T}$.}
\label{thre}
\end{figure}

% \emph{Example:} In Fig. \ref{thre}, a threshold graph $G$ associated with the construction sequence $T^{G}=(0,1,0,1,0,0,1)$ and its corresponding cotree $\mathcal{T}$ are shown.
 \emph{Example:} 
In Fig. \ref{thre}, a threshold graph $G$ and its corresponding cotree $\mathcal{T}$ are shown.
%

  %Now, let us start with one isolated node as the initial graph, and in each step, connect an isolated node to the former graph through the join or the union operation. The obtained graph is referred to as a \emph{threshold graph}, which is a spacial type of a cograph. One can associate a binary \emph{construction sequence} $T^G\in\{0,1\}^n$ to a threshold graph $G$ of size $n$, where $T^G(1)=0$, and for $1<i\leq n$, $T^G(i)=0$ (resp., $T^G(i)=1$) if the node $i$ is added to the former graph by the union (resp. join) operation \cite{hagberg2006designing}. In fact, any threshold graph $G$ can be uniquely determined by its construction sequence $T^G$. In this paper, we assume that all threshold graphs are described and given by their associated construction sequences. We this notation, we have also the next theorem that connects the spectrum of a threshold graph to its degree sequence. 
 %

\subsection{Graph Symmetry}
The symmetries in a graph can be described by its automorphism group~\cite{rahmani2009controllability, zhang2014upper,aguilar2017almost}. Let $\sigma:V\rightarrow V$ be a permutation of the nodes of a graph $G=(V,E)$, where $\{u,v\}\in E$ if and only if $\{\sigma(u),\sigma(v)\}\in E$. Then such a $\sigma$ is referred to as an automorphism of $G$. A trivial automorphism of $G$ is a permutation fixing all nodes. 

\color{black}
\subsection{Problem Formulation}

In this paper, we consider a linear time-invariant (LTI) system defined on a connected cograph $G=(V,E)$ with the Laplacian dynamics described as,
\begin{equation}
\dot{x}=Ax+Bu,
\label{e1}
\end{equation} 
where $A=-L(G)$, and $L(G)\in \mathbb{R}^{n\times n}$ is the Laplacian matrix associated with $G$. Moreover, $x=[x_1,\ldots,x_n ]^T$ is the aggregated vector of states of the nodes\footnote{For notational simplicity, the state of each node is assumed to be a scalar; extension to multi-dimensional case is facilitated using Kronecker products.}, and $u=[u_1,\ldots,u_m ]^T$ is the vector of inputs. Also, $B\in \mathbb{R}^{n\times m}$ is the input matrix whose nonzero entries determine the nodes where the input signals are directly injected. Here, we assume that any input signal can be injected into only one node, referred to as a \emph{control node}. Thus, $B$ assumes the form,
\begin{equation}
B=[e_{j_1},\ldots,e_{j_m}],
\label{B}
\end{equation}
 where $j_i\in\{1,\ldots,n\}$, for $1\leq i\leq m$. We refer to $V_C=\{j_1,\ldots,j_m\}$ as the set of control nodes.
 
 %is assumed to be  a threshold graph or a general cograph, and the controllability of the network is investigated. 
The controllability of an LTI system  captures the ability of an input to steer the states of the system from an arbitrary initial value to any final value within a finite time.  \color{black} In this paper, we aim to provide graph-theoretic controllability conditions for the network  described in (\ref{e1}) and find the minimum number of control nodes from which the network is controllable. 
The celebrated Popov-Belevitch-Hautus (PBH) test has proved
to be instrumental in bridging controllability analysis for networks
to graph-theoretic constructs.
  
 \begin{pro}[\cite{sontag2013mathematical}]
 A system with dynamics (\ref{e1}) (or the pair $(A,B)$)  is controllable if and only if for any (left) eigenvector $\nu$ of  $A$, we have $\nu^TB\neq 0$. 
 \label{prop1}
 \end{pro}
 
% The PBH test can be stated in another equivalent way: a system with dynamics (\ref{e1}) is controllable if and only if for every eigenvalue $\lambda$ of $A$, the matrix $[\lambda I_n-A, B]\in \mathbb{R}^{n\times (n+m)}$ is full rank.
 
Note that if we would like to select a set of control nodes for a network of size $n$ based on the PBH test, we can perform a brute-force verification of the required controllability condition for exponentially many combinations, a computationally impractical endeavor for large-scale networks. Thereby, in this paper, we aim to characterize controllability conditions that can be efficiently inferred from the network topology, specifically for cographs. The key ingredient for such characterizations is the intricate connection between the modularity of cographs and their modal properties.

The following proposition is one of the pertinent results, providing a necessary condition for controllability of a general network based on the existing symmetries in its structure. 

\begin{pro}[\cite{rahmani2009controllability,chapman2015state}]
A network with dynamics (\ref{e1}) defined on a graph $G$ is uncontrollable if $G$ has a nontrivial automorphism that fixes the control nodes.
\label{sym}
\end{pro}

%we show under what operations on network its s-controllability is preserved.

 \color{black}
 
 \section{Cographs Eigenvalues and Eigenvectors} \label{cograph-spectra}

In this section, we investigate the spectrum and an associated modal matrix of a cograph $G$. Indeed, in order to use the PBH test for controllability analysis of a cograph, we first characterize the eigenspace associated with any  eigenvalue of $G$. Then, we find some conditions on the multiplicity of any eigenvalue, which denotes the dimension of the corresponding eigenspace. These results turn out to be instrumental to provide the controllability conditions for a cograph. 

\subsection{Computing Eigenvalues and Eigenvectors of a Cograph} 
In this part, we compute the eigenvalues and linearly independent eigenvectors of a cograph. Indeed, to any internal node $v$ of a cotree, we associate a set of eigenvectors and an eigenvalue whose multiplicity is one less than the number of children of $v$. 
 The next result is an extension of Theorem 2.1 of \cite{merris1998laplacian}.

% provides the nontrivial spectrum and the modal matrix of the join and the union of $p$ graphs $G_1, G_2, \ldots,G_p$ based on the nontrivial spectrum and the modal matrix of each.
%
First, let $\Lambda_i=\Lambda(G_{i})$ and $\mathcal{V}_{i}=\mathcal{V}(G_i)$ be respectively, the nontrivial spectrum and the associated nontrivial modal matrix of the graph $G_i$, $i=1,\ldots,p$ (note that for a 
graph $G_i$ of size one, $\Lambda(G_i)=\emptyset$ and $\mathcal{V}(G_i)=\emptyset$
). Moreover, assume that the nodes of $G_i$ are indexed prior to nodes of $G_j$, for $1\leq i<j\leq p$.% With a slight abuse of notation in the next theorem, for a vector $v\in \mathbb{R}^n$ and a scalar $m\in \mathbb{R}$, we let $v+m=v+m\textbf{1}_{n}$. 
\begin{thm}%[\cite{merris1998laplacian}]
Consider the graphs $G_1,G_2,\ldots,G_p$ of respectively, size $n_1,n_2,\ldots,n_p$, and let $n=\sum_{i=1}^p n_i$. Then,  
%\begin{equation*}L(G_1+G_2)=\begin{bmatrix}L(G_1) & \textbf{0}_{n_1\times n_2}\\ \textbf{0}_{n_2\times n_1} & L(G_2)\end{bmatrix},\end{equation*} 
 %\begin{equation*}L(G_1*G_2)=\begin{bmatrix}L(G_1)-n_2I_{n_1} & \textbf{0}_{n_1\times n_2}\\ \textbf{0}_{n_2\times n_1} & L(G_2)-n_1I_{n_2}\end{bmatrix},\end{equation*}
\begin{itemize} 
\item $\Lambda(G_1+G_2+\ldots+G_p)=(\Lambda_1,\Lambda_2, \ldots,\Lambda_p, 0^{(p-1)})$,\\
\item $\Lambda(G_1*G_2*\ldots*G_p)=\\
(\Lambda_1+n-n_1, \Lambda_2+n-n_2,\ldots, \Lambda_p+n-n_p, n^{(p-1)})$,\\
\item $\mathcal{V}(G_1+G_2+\ldots+G_p)=\mathcal{V}(G_1*G_2*\ldots*G_p)=\\
 \begin{bmatrix} 
\mathcal{V}_1 & \textbf{0}_{n_1} & \ldots & \textbf{0}_{n_1} & n_2\textbf{1}_{n_1} & \ldots & n_p\textbf{1}_{n_1}\\
\textbf{0}_{n_2} & \mathcal{V}_2 & \ldots & \textbf{0}_{n_2} & -n_1\textbf{1}_{n_2} & \ldots & n_p\textbf{1}_{n_2}\\
\vdots & \vdots & \ddots & \vdots & \vdots & \ddots & \vdots\\
\textbf{0}_{n_p} & \textbf{0}_{n_p} & \ldots & \mathcal{V}_p & \textbf{0}_{n_p} & \ldots & -(\sum_{i=1}^{p-1} n_i )\textbf{1}_{n_p}
\end{bmatrix}$.
\end{itemize}
%\begin{equation*}\bar{V}(G_1*G_2)=\bar{V}(G_1+G_2).\end{equation*}
\label{th3}
\end{thm}

\emph{Proof:} The proof is based on an inductive argument. For two graphs $G_1$ and $G_2$, one has ${\Lambda}(G_1+G_2)=(\Lambda_1, \Lambda_2,0)$ and ${\Lambda}(G_1*G_2)=(\Lambda_1, \Lambda_2,n_1+n_2)$ \cite{merris1998laplacian}. Moreover, 
\begin{equation}\mathcal{V}(G_1+G_2)=\mathcal{V}(G_1*G_2)=\begin{bmatrix} \mathcal{V}_1 & \textbf{0}_{n_1} & n_2\textbf{1}_{n_1}\\ \textbf{0}_{n_2} & \mathcal{V}_2 & -n_1\textbf{1}_{n_2}\end{bmatrix}.
\label{eq3}
\end{equation}
Thus the statement of the theorem holds for $p=2$. Now, assume that for $p=k$,  the statement of the theorem is valid. We want to prove the claim for $p=k+1$. Consider $G_1+G_2+\ldots+G_{k}+G_{k+1}$ (respectively, $G_1*G_2*\ldots*G_{k}*G_{k+1}$) as $G'+G_{k+1}$ (respectively, $G'*G_{k+1}$), where $G'=G_1+G_2+\ldots+G_{k}$ (respectively, $G'=G_1*G_2*\ldots*G_{k}$). Thereby, using (\ref{eq3}), the statement of the theorem is valid for $p=k+1$.
\carre
%%%%%%%%

Before characterizing the nontrivial spectrum and an associated nontrivial modal matrix of a cograph, let us introduce more notation. For an internal node $v$ in a cotree $\mathcal{T}=(V^{\mathcal{T}},E^{\mathcal{T}})$, we recall that $\mathcal{C}(v)$ with $c(v)=|\mathcal{C}(v)|$, and $\mathcal{L}(v)$ with $l(v)=|\mathcal{L}(v)|$ are respectively, the set of children and leaves descending from $v$. Let $c(v)=k$, and $\mathcal{C}(v)=\{v_1,\ldots,v_{k}\}$. Note that $k>1$. Now, define $ \lambda_{\mathrm{new}}(v)=\mathrm{lab}(v)\times l(v)$, which is referred to as the new eigenvalue of the internal node $v$. Then, if $\mathrm{lab}(v)=0$, $ \lambda_{\mathrm{new}}(v)=0$, and  if $\mathrm{lab}(v)=1$,  $ \lambda_{\mathrm{new}}(v)=l(v)$. Now, let $n_i=l(v_i)$, $1\leq i\leq k$, and consider the matrix $M\in \mathbb{R}^{l(v)\times (k-1)}$  as %Consider a matrix $M\in \mathbb{R}^{l(v)\times c(v)-1}$ whose $i$th column, $1\leq i\leq c(v)-1$, is defined as $M_{:,i}=[n_{i+1}\textbf{1}_{\sum_{k=1}^i n_k}^T, -(\sum_{k=1}^i n_k)\textbf{1}_{n_{i+1}}^T, \textbf{0} ]^T$. Let us define $\mathcal{V}_{\mathrm{new}}(v)=M$. Then, it can be written as
\begin{equation*}
\scriptsize
M=\begin{bmatrix}
n_2 \textbf{1}_{n_1} & \ldots & n_{k-1}\textbf{1}_{n_1} & n_{k}\textbf{1}_{n_1} \\
-n_1\textbf{1}_{n_2} & \ldots & n_{k-1}\textbf{1}_{n_2} & n_{k}\textbf{1}_{n_2} \\
\vdots & \ddots & \vdots & \vdots\\
\textbf{0}_{n_{k-1}} & \ldots & -\sum_{j=1}^{k-2} n_j \textbf{1}_{n_{k-1}} & n_{k}\textbf{1}_{n_{k-1}}\\
\textbf{0}_{n_{k}} & \ldots & \textbf{0}_{n_{k}}  & -\sum_{j=1}^{k-1} n_j \textbf{1}_{n_{k}}\\
\end{bmatrix}.
\end{equation*}
Let us define $\mathcal{V}_{\mathrm{new}}(v)= M$, which we refer to as a new modal matrix of the internal node $v$.

Now, consider two internal nodes $v,w\in V^{\mathcal{T}}$, where $w$ is an ancestor of $v$. Let $\mathcal{P}_v^w=(u_0,u_1,\ldots,u_p)$, where $u_0=w$ and $u_p=v$. Then, for an eigenvalue $\lambda\in \Lambda(G_{(v)})$,  if $v=w$, the updated eigenvalue of $v$ at $w$ is defined as $\lambda_{\mathrm{upd}}^w(v)=\lambda_{\mathrm{new}}(v)$; otherwise, it is defined as
 $$\lambda_{\mathrm{upd}}^w(v)=
\lambda_{\mathrm{new}}(v)+\sum_{i=0}^{p-1} \mathrm{lab}(u_{i})\times(l(u_{i})-l(u_{i+1})).
$$
 In addition,  let $\mathcal{V}'\in\mathbb{R}^{l(w)\times (c(v)-1)}$ be such that $\mathcal{V}'_{\mathcal{L}(v),:}=\mathcal{V}_{\mathrm{new}}(v)$, and $\mathcal{V}'_{\mathcal{L}(w)\setminus \mathcal{L}(v),:}=\textbf{0}_{(l(w)-l(v))\times (c(v)-1)}$. In other words, $\mathcal{V}'$ is a  $l(w)\times (c(v)-1)$ matrix whose rows corresponding to the indices of leaves of $G_{(v)}$ constitute the matrix  $\mathcal{V}_{\mathrm{new}}(v)$, while the rest of its rows are the zero vectors.  Let us define the updated modal matrix of $v$ at $w$ as $\mathcal{V}_{\mathrm{upd}}^w(v)=\mathcal{V}'$. It is obvious that $\mathcal{V}_{\mathrm{upd}}^v(v)=\mathcal{V}_{\mathrm{new}}(v)$.
 
 \begin{thm}
Consider a cograph $G=(V,E)$ with the associated cotree $\mathcal{T}=(V^{\mathcal{T}},E^{\mathcal{T}} )$ and the root $r$. Let $h=|V^{\mathcal{T}}|-|V|$ be the number of internal nodes of $\mathcal{T}$, and let $v^1, \ldots, v^h$ be its internal nodes. For $1\leq i\leq h$, let $\lambda_i=\lambda_{\mathrm{upd}}^r(v^i)$, and $\mathcal{V}^{(i)}=\mathcal{V}_{\mathrm{upd}}^r(v^i)$. Then,\\
%\begin{itemize}
\begin{equation} \Lambda(G)=(\lambda_1^{(c(v^1)-1)},\ldots,\lambda_h^{(c(v^h)-1)}),\label{spect}\end{equation}
 \begin{equation}\mathcal{V}(G)=[\mathcal{V}^{(1)},\ldots,\mathcal{V}^{(h)}].\label{modal}\end{equation}
%\end{itemize} 
\label{thmm4}
\end{thm}

\emph{Proof:} The proof follows by a strong induction on $h$. %The theorem is proved using a strong inductive argument.
 First, we show that the result holds for $h=1$. Let $v$ be the single internal node with $c(v)=l(v)=n$, and note that $r=v$.  Then, the children of $v$ are all graphs of size one. Thus, from Theorem \ref{th3}, if $\mathrm{lab}(v)=0$, $\Lambda(G)=(0^{(n-1)})$, and otherwise $\Lambda(G)=(n^{(n-1)})$. Hence, $\lambda_{\mathrm{upd}}^r(v)=\lambda_{\mathrm{new}}^r(v)=\mathrm{lab}(v)\times l(v)$, and for $\lambda=\lambda_{\mathrm{upd}}^r(v)$, one can write $\Lambda(G)=(\lambda^{(c(v)-1)})$. Moreover, Theorem \ref{th3} implies that 
\begin{equation}\scriptsize
\mathcal{V}(G)=\begin{bmatrix}
1 & \ldots & 1 & 1\\
-1 & \ldots & 1 & 1\\
\vdots & \ddots & -(n-2) & 1\\
0 & \ldots & 0 & -(n-1)
\end{bmatrix}
\end{equation}
Consequently, one can write $\mathcal{V}(G)=\mathcal{V}_{\mathrm{upd}}^r(v)$; the result is thus valid for $h=1$. Now, assuming that the result holds for all $h\leq k$, we want to prove that it holds for $h=k+1$. Let $r=v^{k+1}$, and $c(r)=p$. Further, let $\mathcal{C}(r)=\{u_1,\ldots,u_p\}$. Let us index the leaves of $\mathcal{T}$ in a way that for $1\leq i<j\leq p$, leaves of $\mathcal{T}_{(u_i)}$ are indexed prior to the leaves of $\mathcal{T}_{(u_j)}$. Since the number of internal nodes of every $\mathcal{T}_{(u_i)}$, $1\leq i\leq p $, is less than $k+1$, by our inductive hypothesis, we know that  $\Lambda(G_{(u_i)})$ is a sequence of $\lambda_{\mathrm{upd}}^{u_i}(w)$ with the multiplicity $c(w)-1$, where $w$ is an internal node of $\mathcal{T}_{(u_i)}$. Then, from Theorem \ref{th3}, $\Lambda(G)$ includes a sequence of $\lambda_{\mathrm{upd}}^{u_i}(w)+\mathrm{lab}(r)(l(r)-l(u_i))=\lambda_{\mathrm{upd}}^{r}(w)$ with multiplicity $c(w)-1$ for every $w$ which is an internal node of one of $\mathcal{T}_{(u_i)}$, $1\leq i\leq p$. Moreover, $\Lambda(G)$ includes the eigenvalue $\lambda_{\mathrm{new}}(r)=\lambda_{\mathrm{upd}}^{r}(v^{k+1})$ with the multiplicity $c(r)-1$. Thus, the result is valid for the nontrivial spectrum of $G$, when $h=k+1$. Using a similar argument, based on the inductive assumption, $\mathcal{V}(G_{(u_i)})$, $1\leq i\leq p$, is a sequence of $\mathcal{V}_{\mathrm{upd}}^{u_i}(w)$, where  $w$ is an internal node of $\mathcal{T}_{(u_i)}$. In addition, Theorem \ref{th3} implies that $\mathcal{V}(G)$ includes $\mathcal{V}_{\mathrm{upd}}^{r}(w)$, for every $w$ that is an internal node of one of $\mathcal{T}_{(u_i)}$'s, $1\leq i\leq p$. Also, $\mathcal{V}(G)$ includes $\mathcal{V}_{\mathrm{new}}(r)=\mathcal{V}_{\mathrm{upd}}^{r}(v^{k+1})$, and thereby, the result is valid for $h=k+1$. Thus, the assertion  holds for any cograph.\carre   %A similar argument can also be made for the nontrivial modal matrix of $G$, and then  the result of the theorem can be verified.   \carre 

Using Theorem \ref{thmm4}, one can also find a relationship between the number of leaves of a rooted tree and the number of children of its internal nodes.

\begin{cor}
Let $n$ be the number of leaves of a rooted tree, and $v^1,\ldots, v^h$ be its internal nodes. Then, $n-1=\sum_{i=1}^h (c(v^i)-1)$.
\label{cor4}
\end{cor}

\emph{Proof:} We have $|\Lambda(G)|=n-1$. Moreover, from equation (\ref{spect}), $|\Lambda(G)|=\sum_{i=1}^h (c(v^i)-1)$,  completing the proof.
\carre

Based on Corollary \ref{cor4}, we can also state the next result for a cograph. This result was initially stated in \cite{corneil1981complement}.
% through an almost intuitive proof. 
\begin{pro}
Any cograph $G=(V,E)$, where $|V|>1$, has at least a pair of siblings.
\label{prop2}
\end{pro}

\emph{Proof:} Consider the cotree $\mathcal{T}$ associated with $G$, and let $n$ and $h$ be respectively, the number of leaves and internal nodes of $\mathcal{T}$. Assume that no two nodes in $G$ are siblings. Then, every internal node of $\mathcal{T}$ has at most one child which is a leaf. This implies that $n\leq h$. In addition, for every internal node $v^i$, $1\leq i\leq h$, we have $c(v^i)\geq 2$. Thereby, from Corollary \ref{cor4}, $n-1=\sum_{i=1}^h (c(v^i)-1)\geq h$, contradicting $n\leq h$.  
\carre

\emph{Example:} Considering the cotree in Fig. \ref{cograph} (b), one can find that $\Lambda(G)=(7^{(1)},5^{(1)},8^{(2)},4^{(1)},7^{(2)})$. For example, for internal node $v^5$, we have $\lambda_{\mathrm{new}}(v^5)=l(v^5)=3$. Also, one has $\lambda_5=\lambda^r_{\mathrm{upd}}(v^5)=3+0\times(l(v^4)-l(v^5))+1\times(l(v^3)-l(v^4))=7$, and $c(v^5)-1=2$. Moreover, $\scriptsize\mathcal{V}^T_{\mathrm{new}}(v^5)=\begin{bmatrix}
1 & -1 & 0\\
1 & 1& -2
\end{bmatrix}$. As a result,  a modal matrix associated with $\Lambda(G)$ is obtained as,
$$\scriptsize\mathcal{V}(G)=\begin{bmatrix}
 1 & 1 & 1 & 4 & 0 & 0 & 0 \\
 -1 & 1 & 1 & 4 & 0 & 0 & 0 \\
  0 & -2 & 1 & 4 & 0 & 0 & 0 \\
  0 & 0   & -3 & 4 & 0 & 0 & 0 \\
  0 & 0   & 0   & -4 & 3 & 0 &0  \\
  0 & 0   & 0   & -4 & -1 & 1 & 1 \\
  0 &  0  & 0   & -4 & -1 & -1 & 1 \\
 0  & 0   & 0   & -4 & -1 & 0 & -2 
\end{bmatrix}.$$
\color{black}
\subsection{Conditions on Multiplicity of an Eigenvalue and the Structure of the Modal Matrix}

In this part, we derive conditions under which the eigenvalues associated with two internal nodes of a cotree are distinct. Also, we will see how to choose some rows of a modal matrix associated with an eigenvalue of a cograph such that the resulting matrix is invertible.
  
As mentioned before, one can characterize the nontrivial spectrum of a cograph by using (\ref{spect}) in Theorem \ref{thmm4}. However, note that for some $1\leq i<j\leq h$, we may have $\lambda_i=\lambda_j$. The next result identifies conditions under which the updated eigenvalues of two internal nodes at the root are distinct.

 \begin{lem}
 Let $r$ be the root and $v,w$ be two internal nodes of a cotree. If  $v$ is an ancestor of $w$, then $\lambda_{\mathrm{upd}}^r(v)\neq \lambda_{\mathrm{upd}}^r(w)$.
 \label{lemm1}
 \end{lem}
 
  \emph{Proof:} In a cotree, there is a unique path  from $r$ to $v$ ($\mathcal{P}_{v}^{r}$), and since $v$ is an ancestor of $w$, there is also a unique path  from $v$ to $w$ ($\mathcal{P}_{w}^{v}$). Let $\mathcal{P}_{w}^{v}=(u_0,u_1,\ldots,u_k)$, where $u_0=v$, and $u_k=w$. Now, let us first prove that $\lambda_{\mathrm{upd}}^{u_0}(u_k)\neq \lambda_{\mathrm{new}}(u_0)$. For $0\leq i\leq k$, let $n_i=l(u_i)$, and note that the number of leaves of an internal node is greater than the number of leaves of any of its children, that is, $n_{i}>n_{i+1}$, $0\leq i\leq k-1$. Moreover, in a cotree, the label of nodes of a path  alternates between 0 and 1. Then, if $\mathrm{lab}(u_i)=0$, $\mathrm{lab}(u_{i+1})=1$, and vice versa. Without loss of generality, assume that $k$ is even, say $k=2r$ for some integer $r$. Then, if $\mathrm{lab}(w)=1$,  by the definition, we have $\lambda_{\mathrm{upd}}^{v}(w)=n_k+\sum _{i=0}^{r-1} (n_{2i}-n_{2i+1})$. Moreover, if $\mathrm{lab}(w)=0$,  $\lambda_{\mathrm{upd}}^{v}(w)=\sum _{i=0}^{r-1} (n_{2i+1}-n_{2i+2})$. Hence, since $n_0=n_k+\sum _{i=0}^{r-1} (n_{2i}-n_{2i+1})+\sum _{i=0}^{r-1} (n_{2i+1}-n_{2i+2})$,  it follows that $0<\lambda_{\mathrm{upd}}^{v}(w)< n_0 $. In addition, note that $\lambda_{\mathrm{new}}(v)$ is either 0 or $n_0$. Thus, $\lambda_{\mathrm{upd}}^{v}(w)\neq \lambda_{\mathrm{new}}(v)$. Now, let $\mathcal{P}_{v}^{r}=(z_0,\ldots, z_s)$, where $z_0=r$ and $z_s=v$. Define $\lambda_{\mathrm{add}}=\sum_{i=0}^{s-1} \mathrm{lab}(z_{i})(l(z_{i})-l(z_{i+1}))$. Thereby, we have $\lambda_{\mathrm{upd}}^r(v)=\lambda_{\mathrm{new}}(v)+\lambda_{\mathrm{add}}$, and $\lambda_{\mathrm{upd}}^r(w)=\lambda_{\mathrm{upd}}^v(w)+\lambda_{\mathrm{add}}$. Thus, $\lambda_{\mathrm{upd}}^r(v)\neq \lambda_{\mathrm{upd}}^r(w)$, and the proof is complete.  
 \carre 
 
 Note that from Lemma \ref{lemm1}, the updated eigenvalue of two internal nodes $v,w$ at the root of a cotree may be the same only in the case that none of these nodes is the ancestor of the other one. We now show that in this case,  the index sets of leaves of $v$ and $w$ have an empty intersection.
 \begin{pro}
 Consider a  cotree $\mathcal{T}=(V^{\mathcal{T}},E^{\mathcal{T}} )$. For two nodes $v,w\in V^{\mathcal{T}}$, %and without loss of generality, assume that $v$ is not an ancestor of $w$. 
  $\mathcal{L}(v)\cap \mathcal{L}(w)\neq \emptyset$ if and only if either $w$ is an ancestor of $v$, or $v$ is an ancestor of $w$. 
 \label{prop3}
 \end{pro}
 
 \emph{Proof:} The proof follows by a contradiction. Assume that neither is $w$ an ancestor of $v$, nor is $v$  an ancestor of $w$. Then the lowest common ancestor of $v$ and $w$ is some node, say $z$, where $z\neq v$ and $z\neq w$.  Moreover, we assume that there is some leaf $i\in V^\mathcal{T}$ such that $i\in\mathcal{L}(v)\cap \mathcal{L}(w)$. Since $\mathcal{T}$ is a rooted tree, there should be a unique path from the root $r$ to  the leaf $i$. However, one can find two paths $\mathcal{P}_1=(\mathcal{P}_z^r,\mathcal{P}_v^z,\mathcal{P}_i^v)$ and $\mathcal{P}_2=(\mathcal{P}_z^r,\mathcal{P}_w^z, \mathcal{P}_i^w)$ that are both directed from $r$ to $i$, establishing a contradiction.   
 \carre
 
\emph{Example:} Considering the cotree in Fig. \ref{cograph} (b), we have $\lambda_{\mathrm{upd}}^r(v^1)=\lambda_{\mathrm{upd}}^r(v^5)=7$. Also, $\mathcal{L}(v^1)=\{1,2\}$, and $\mathcal{L}(v^5)=\{6,7,8\}$; thus, $\mathcal{L}(v^1)\cap \mathcal{L}(v^5)=\emptyset$.  

\color{black}
 Now for two internal nodes $v,w$ in a cotree $\mathcal{T} $ associated with the  cograph $G$, let $\lambda_{\mathrm{upd}}^r(v)=\lambda_{\mathrm{upd}}^r(w)=\lambda$. Then for a full rank matrix $\bar{\mathcal{V}}\in \mathbb{R}^{n\times (c(v)+c(w)-2)}$ defined as $\bar{\mathcal{V}}=[\mathcal{V}_{\mathrm{upd}}^r(v),\mathcal{V}_{\mathrm{upd}}^r(w)]$, we have $L(G)\bar{\mathcal{V}}=\lambda \bar{\mathcal{V}}$. From Lemma \ref{lemm1} and Proposition \ref{prop3}, one can then conclude that $\bar{\mathcal{V}}_{\mathcal{L}(v),:}=[\mathcal{V}_{\mathrm{new}}(v),\textbf{0}_{\mathcal{L}(v)\times (c(w)-1)}]$, and $\bar{\mathcal{V}}_{\mathcal{L}(w),:}=[\textbf{0}_{\mathcal{L}(w)\times (c(v)-1)},\mathcal{V}_{\mathrm{new}}(w)]$; while the other rows of $\bar{\mathcal{V}}$ are zero. %$\mathcal{V}_{\{1,\ldots,n\}\setminus (\mathcal{L}(v)\cup\mathcal{L}(w)),:}=\textbf{0}_{}$ 
Let $\mathcal{S}\subset\{1,\ldots,n\}$, and assume that for some $\mathcal{S}_1\subseteq \mathcal{L}(v)$ and $\mathcal{S}_2\subseteq \mathcal{L}(w)$, one has $\mathcal{S}=\mathcal{S}_1\cup \mathcal{S}_2$. Then, it follows that a submatrix of $\bar{\mathcal{V}}$ with rows chosen from indices of $\mathcal{S}$,  denoted by $\bar{\mathcal{V}}_{\mathcal{S},:}$, is nonsingular if and only if $(\mathcal{V}_{\mathrm{upd}}^r(v))_{\mathcal{S}_1,:}$ and $(\mathcal{V}_{\mathrm{upd}}^r(w))_{\mathcal{S}_2,:}$ are both nonsingular. In this case, we present conditions under which for an internal node $v$,  $(\mathcal{V}_{\mathrm{upd}}^r(v))_{\mathcal{S},:}$ is invertible.

\emph{\bf Procedure I:} Given an internal node $v$ in a cotree, first choose a subset $\mathcal{S}'$ of children of $v$, where $|\mathcal{S}'|=c(v)-1$. Then, select a leaf from $\mathcal{L}(u)$ for any $u\in \mathcal{S}'$. Let $\mathcal{S}$ be the set of selected leaves. 

\begin{lem} %Consider an internal node $v$ in a cotree, where $k=c(v)$ and $\mathcal{C}(v)=\{v_1,\ldots,v_{k}\}$. Define $\mathcal{V}=\mathcal{V}_{\mathrm{upd}}^r(v)$, and let $s_1=\{j_1,\ldots,j_{k-1}\}$, where $1\ldots j_i\ldots n$, for $1\leq i\leq k-1$. Also, let $s_2=\{v_{j_1},\ldots,v_{j_{k-1}}\}, where $v_{j_i}\in\mathcal{C}(v)$.  Then, $\mathcal{V}_{s,:}$ is nonsingular if and only if any $j_i\in s_1$ corresponds to one leaf of  a $v$ 
%Consider an internal node $v$ in a cotree, where $k=c(v)$ and $\mathcal{C}(v)=\{v_1,\ldots,v_{k}\}$.  Define $\mathcal{V}=\mathcal{V}_{\mathrm{upd}}^r(v)$, where $\mathcal{V}\in \mathbb{R}^{n\time (k-1)}$. Let $P\in \mathbb{R}^{(k-1)\times(k-1)}$ be a submatrix of $\mathcal{V}$. Then, $P$ is nonsingular if and only if any of its rows corresponds to one leaf of a  Define $\mathcal{V}=\mathcal{V}_{\mathrm{upd}}^r(v)$, and let $s_1=\{j_1,\ldots,j_{k-1}\}$, where $1\ldots j_i\ldots n$, for $1\leq i\leq k-1$. Also, let $s_2=\{v_{j_1},\ldots,v_{j_{k-1}}\}, where $v_{j_i}\in\mathcal{C}(v)$.  Then, $\mathcal{V}_{s,:}$ is nonsingular if and only if any $j_i\in s_1$ corresponds to one leaf of  a $v$
Let $\mathcal{V}=\mathcal{V}_{\mathrm{upd}}^r(v)$, where $v $ is an internal node of a cotree. Then  $\mathcal{V}_{\mathcal{S},:}$ is nonsingular if and only if $\mathcal{S}$ is chosen according to Procedure I.
\label{lem5}
\end{lem}

\emph{Proof:} Let $k=c(v)$, and assume that $\mathcal{V}_{\mathcal{S},:}\in \mathbb{R}^{(k-1)\times (k-1)}$ is nonsingular. Then, it follows that $\mathcal{S}\subset \mathcal{L}(v)$. Otherwise, $\mathcal{V}_{\mathcal{S},:}$ has some zero rows, establishing a contradiction. %Moreover,  none of two rows of $\mathcal{V}_{s,:}$ should be the same. Then, by considering the definition of $\mathcal{V}_{\mathrm{upd}}^r(v)$, we should choose at most one leaf from any child $u$ of $v$. Now, consider the matrix $J$ defined as
Moreover, for any child $u$ of $v$, all the rows of $\mathcal{V}$ corresponding to the leaves of $u$ are the same. Then, we should choose at most one leaf from any child of $v$. Let $\mathcal{C}(v)=\{v_1,\ldots,v_k\}$, and $n_i=l(v_i)>0$, $1\leq i\leq k$. Consider the matrix $M^{(k)}\in \mathbb{R}^{k\times (k-1)}$ defined as
\begin{equation*}
\scriptsize
M^{(k)}= \begin{bmatrix} 
n_2 & \ldots & n_{k-1} & n_k\\
-n_1 & \ldots & n_{k-1} & n_k\\
\vdots & \ddots &\vdots & \vdots\\
0 & \ldots & -\sum_{i=1}^{k-2} n_i & n_k\\
0 & \ldots & 0 & -\sum_{i=1}^{k-1} n_i 
\end{bmatrix},
\end{equation*}
whose arbitrary row corresponds to one leaf of a child of $v$. 
It now suffices to show that by choosing any $k-1$ rows of $M^{(k)}$, a nonsingular $(k-1)\times(k-1)$ matrix is obtained.
The proof is based on an induction on $k$. Let $k=2$. Then, $M^{(2)}=[n_2, -n_1]^T$, and any of its $1\times 1$ submatrices is nonzero and nonsingular. Now, assume that for $k=h$, and for all $n_i>0$, $1\leq i\leq h$, any $(h-1)\times (h-1)$ submatrix of $M^{(k)}$ is  nonsingular. Based on this assumption,  we claim  that for $k=h+1$, any $h$ rows of $M^{(k)}$ are linearly independent. Let $R$ with $|R|=h$ be the indices of the rows chosen. First, assume that $R=\{i_1,\ldots, i_{h-1}, h+1\}$, where $1\leq i_j\leq h$, for $1\leq j\leq h-1 $. 
%$$ where $r_i\in \mathbb{R}^{1\times(h-1)}$, and $[r_1^T,r_2^T,\ldots,r_h^T]^T=M^{(h)}$. Then, based on the inductive assumption, every $h-1$ rows  of $M^{(h)}$ are linearly independent.  
Moreover, let $R'=\{i_1,\ldots, i_{h-1}\}$. Then we can write $$\scriptsize M^{(h+1)}_{R,:}=\begin{bmatrix}
M^{(h)}_{R',:} & n_{k+1}\textbf{1}_{h-1}\\
\textbf{0}_{1\times (h-1)} & -\sum_{i=1}^{h} n_i\end{bmatrix}.$$
If $M^{(h+1)}_{R,:}$ is singular, there is a nonzero $\alpha \in \mathbb{R}^{1\times h}$ such that $\alpha M^{(h+1)}_{R,:} =0$. Let $\alpha=[\beta,\alpha_h]$, where $\beta\in \mathbb{R}^{h-1}$, and $\alpha_h\in \mathbb{R}$. Then, we have $\beta M^{(h)}_{R',:} =0$, and since from the inductive assumption,  $M^{(h)}_{R',:} $ is nonsingular, we can conclude that $\beta=0$. Moreover, $-\alpha_{h}\sum_{i=1}^{h}n_i=0$, which leads to $\alpha_h=0$. Thus, $\alpha=0$, and $M^{(h+1)}_{R,:}$ is nonsingular.  
 Now, let $R=\{1,\ldots,h\}$ and $R'=\{1,\ldots,h-1\}$. By assuming that $\mathrm{det}(M^{(h)}_{R',:} )\neq 0$, we should show that $\mathrm{det}(M^{(h+1)}_{R,:} )\neq 0$. Let $D^*=\mathrm{det}(M^{(h)}_{R',:} )$. One can now verify that $\mathrm{det}(M^{(h+1)}_{R,:} )=\frac{n_{h+1}}{n_{h}} (\sum_{i=1}^{h-1}n_i) D^*+n_{h+1} D^*$ which is nonzero. Hence, $M^{(h+1)}_{R,:}$ is nonsingular, thus completing the proof.      
\carre

\section{Controllability of Cographs} \label{controllability}

In this section, we investigate the controllability of a network with dynamics (\ref{e1}) and the input matrix (\ref{B}), defined on a cograph $G$. First, let us introduce a sibling partition in a cograph.

%\subsection{Sibling Partition}

Consider an undirected graph $G=(V,E)$, and let  $\mathbf{C}_i$, $1\leq i\leq p$, be a nonempty subset of $V$ called a cell. Then, $\pi=\{{\mathbf{C}_1},\ldots,\mathbf{{C}}_p\}$ is a \emph{partition} of $G$ if $\bigcup_{i=1}^p \mathbf{C}_i= V$, and $\mathbf{C}_i\cap \mathbf{C}_j=\emptyset$, for $1\leq i<j\leq p$. Now, let $G$ be a cograph, and let $\pi$ be a partition where any two nodes $i,j\in V$ are siblings if and only if for some $1\leq k\leq p$, $i,j\in \mathbf{C}_k$. Then, we refer to $\pi$ as the \emph{sibling partition} of $G$ and denote it by $\pi_{\mathrm{sib}}(G)$. Note that by this definition, for a cograph $G$, $\pi_{\mathrm{sib}}(G)$ is unique.  Once a cotree associated with a cograph is built, its sibling partition can be found in $\mathcal{O}(n)$. 

\begin{figure}[hbt]
\includegraphics[width=.45\textwidth]{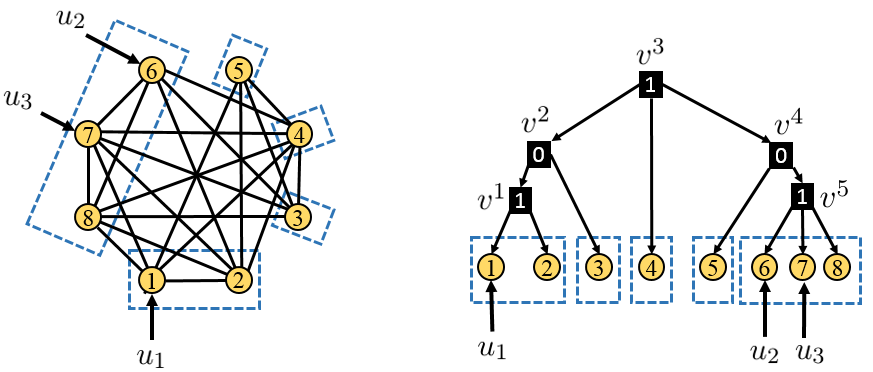}
\centering
\caption{Sibling partition and control nodes in a cograph and its cotree.}
\label{control}
\end{figure}

 \emph{Example:} For the cograph $G$ shown in Fig. \ref{control}, we have $\mathbf{C}_1=\{1,2\}$, $\mathbf{C}_2=\{3\}$, $\mathbf{C}_3=\{4\}$, $\mathbf{C}_4=\{5\}$, and  $\mathbf{C}_5=\{6,7,8\}$\color{black}. One can see that any cell includes all the leaves with the same parent.
% $\pi_{\mathrm{sib}}(G)=\{\{1,2\},\{3\},\{4\}, \{5\},\{6,7,8\}\}$.

\color{black}

\begin{pro} \label{prop-symmetry}
Any permutation that permutes two nodes of a cell in a sibling partition of a cograph $G$ and fixes all other nodes is an automorphism of $G$.
\label{auto}
\end{pro}

\emph{Proof:} Let $\sigma:V\rightarrow V$ be a permutation such that for the two nodes $i,j\in V$, $\sigma(i)=j$ and $\sigma(j)=i$, and for any node $k\notin\{i,j\}$, $\sigma(k)=k$. Then, the edge $\{k,l\}\in E$ is mapped to itself if either $k,l\notin \{i,j\}$ or $k=i$ and $l=j$.  Otherwise, any edge $\{i,k\}\in E$ (respectively, $\{j,k\}\in E$), where $k\notin \{i,j\}$, is mapped to $\{\sigma(i),\sigma(k)\}=\{j,k\}$ (respectively, $\{\sigma(j),\sigma(k)\}=\{i,k\}$). Moreover, since $i$ and $j$ are siblings, $k\in N(i)$ if and only if  $k\in N(j)$. Thus, for any edge $\{l,z\}\in E$, where $l,z\in V$, we have $\{\sigma(l),\sigma(z)\}\in E$, and hence, $\sigma$ is an automorphism of $G$.  
\carre
\begin{rem}
Note that the sibling partition $\pi_{\mathrm{sib}}(G)$ is an equitable partition of a cograph in the sense that all nodes in any cell of $\pi_{\mathrm{sib}}(G)$ have the same neighbors in other cells (for more information about equitable partitions, one may refer to \cite{rahmani2009controllability, zhang2014upper,aguilar2017almost}).  
\end{rem}

\color{black}
The next theorem, which is the main result of this paper, presents a procedure for selecting a minimal set of control nodes in a Laplacian network defined on a cograph. In fact, in order to break the symmetries in every cell $\mathbf{C}_i$ of size $m_i$ in the sibling partition of a cograph, one needs to choose $m_i-1$ control nodes from that cell. Next, we show that  the resulting set of control nodes ensures the network controllability.  
\color{black}
\begin{thm} \label{min-control-cographs}
Consider a network defined on a connected cograph $G$ of size $n>1$ with dynamics (\ref{e1}). Let $\pi_{\mathrm{sib}}(G)=\{ \mathbf{\mathbf{C}}_1,\ldots,\mathbf{C}_p\}$, where $|\mathbf{C}_i|=m_i$, $1\leq i\leq p$. 
 Then, the minimum number of control nodes rendering the network controllable is $n-p$.
 Moreover, a control node set of size $n-p$ should be chosen by selecting $m_i-1$ nodes from any cell $\mathbf{C}_i$ of $\pi_{\mathrm{sib}}(G)$, $1\leq i\leq p$. %$\mbox{min}|V_C|=n-p$.
 \label{thm2}  
\end{thm}

%Note that Theorem \ref{thm2} presents a procedure for selecting a minimal set of control nodes in a cograph.

%We are now ready to prove Theorem \ref{thm2}.

\emph{Proof:}  First, assume that the network is controllable, yet $|V_C|<n-p$. Then, there is a cell $\mathbf{C}_k$, for some $1\leq k\leq p$, at least two nodes of which, say $i$ and $j$, are not in $V_C$. Now, consider a permutation $\sigma$ of $G$ that permutes $i$ and $j$ and fixes all other nodes of the graph (including all control nodes). Based on Proposition \ref{auto}, $\sigma$ is an automorphism of $G$, and thus from Proposition \ref{sym}, the system is uncontrollable, contradicting the assumption.    
\color{black}

%Consider the cotree $\mathcal{T}$ with the root $r$, associated with $G$. Let $v^1,\ldots, v^h$ be the internal nodes of $\mathcal{T}$. Then from Theorem \ref{thmm4}, the nontrivial spectrum and  modal matrix of $G$ are obtained by  (\ref{spect}) and (\ref{modal}). First, assume that the network is controllable from a set of control nodes $V_C$; however, $|V_C|<n-p$. Then $V_C$ does not include any set $S_C$ chosen by Procedure II, as $|S_C|=n-p$. Then from Lemma \ref{lemlast}, there is an internal node $v_i$, $1\leq i\leq h$, for which $B^T\mathcal{V}_{\mathrm{upd}}^r(v_i)$ is not full rank. Accordingly, there is a nonzero vector $\alpha\in \mathbb{R}^n$ such that $B^T\mathcal{V}_{\mathrm{upd}}^r(v_i)\alpha=0$. Note that  $\nu=\mathcal{V}_{\mathrm{upd}}^r(v_i)\alpha$ is a nonzero eigenvector associated with the eigenvalue $\lambda_{\mathrm{upd}}^r(v_i)$. In other words, $G$ has a nonzero eigenvector $\nu$ where $\nu^TB=0$. This however, according to Proposition \ref{prop1}, implies that the network is not controllable, establishing a contradiction. Thus $|V_C|\geq n-p$.

Now,  consider \color{black} the cotree $\mathcal{T}$ with the root $r$, associated with $G$. Let $v^1,\ldots, v^h$ be the internal nodes of $\mathcal{T}$. Then from Theorem \ref{thmm4}, the nontrivial spectrum and  modal matrix of $G$ are obtained by  (\ref{spect}) and (\ref{modal}). Let  $V_C$ be obtained by choosing $m_i-1$ nodes from any cell $\mathbf{C}_i$, $1\leq i\leq p$; however, assume that the network is not controllable. Therefore, from Proposition \ref{prop1}, there should be a nonzero eigenvector $\nu$ associated with the eigenvalue $\lambda$, where $\nu^TB=0$. For $i_j\in\{1,\ldots, h\}$, $1\leq j\leq k$, assume that $\lambda_{i_j}=\lambda$, where $\lambda_{i_j}=\lambda_{\mathrm{upd}}^r(v^{i_j})$.  Now define $\mathcal{V}=[\mathcal{V}^{(i_1)},\ldots,\mathcal{V}^{(i_k)}]$, where for $1\leq j\leq k$, $\mathcal{V}^{(i_j)}=\mathcal{V}_{\mathrm{upd}}^r(v^{i_j})$, and note that $L(G)\mathcal{V}=\lambda \mathcal{V}$. Hence, for some nonzero $\alpha\in \mathbb{R}^n$, one can write $\nu=\mathcal{V}\alpha$. Moreover,  note that from Lemma \ref{lemm1} and Proposition \ref{prop3}, for $1\leq j<l\leq k$, we have $\mathcal{L}(v^{i_j})\cap \mathcal{L}(v^{i_l})=\emptyset$. This simply implies that % for $1\leq i_j\leq k$, $\mathcal{V}_{\mathcal{L}(v_{i_j}),:}=[\textbf{0},\ldots,\textbf{0},\mathcal{V}_{\mathrm{new}}(v_{i_j}), \textbf{0},\ldots,\textbf{0}]$, 
 $B^T\mathcal{V}$ is full rank if and only if $B^T\mathcal{V}^{(i_j)}$, for every $1\leq j\leq k$, is full rank. 

Now, we want to prove that by choosing this $V_C$, $B^T\mathcal{V}^{(i_j)}$ is full rank. Let $v^{i_j}$ be an internal node  of $\mathcal{T}$ with $\mathcal{C}(v^{i_j})=V_{\mathrm{int}}\cup V_{\mathrm{leaf}}$, where $V_{\mathrm{int}}$ (respectively, $V_{\mathrm{leaf}}$) is  the set of children of $v^{i_j}$ that are internal nodes (respectively, leaves) of $\mathcal{T}$. Let $u\in V_{\mathrm{int}}$. Then, from Proposition \ref{prop2}, the cograph $G_{(u)}$  has at least two nodes which are siblings. Accordingly, by choosing this set of control nodes $V_C$, there is a leaf $w\in \mathcal{L}(u)$ such that $w\in V_C$.  Moreover,  %, while $V_{\mathrm{leaf}}$ is the set of children of $v$ that are leaves of $\mathcal{T}$. %$u_i$, $1\leq i\leq p$, is an internal node of $\mathcal{T}$, while $l_j$, $1\leq j\leq q$ is a leaf of $\mathcal{T}$. If   
 if $V_{\mathrm{leaf}}\neq \emptyset$, it includes the leaves of $\mathcal{T}$ with the same parent $v^{i_j}$, and thereby, for some $1\leq k\leq p$, $V_{\mathrm{leaf}}=\mathbf{C}_k$, where $\mathbf{C}_k$ is a cell of the sibling partition, and $|V_{\mathrm{leaf}}|=m_k$. Hence, for the internal node $v^{i_j}$, $V_C$ includes one leaf of any internal node $u\in V_{\mathrm{int}}$. Moreover, it includes $m_k-1$ nodes of $V_{\mathrm{leaf}}$. Therefore, for every internal node $v^{i_j}$, $V_C$ includes an associated  set $\mathcal{S}$ chosen by Procedure I, and thus Lemma \ref{lem5} implies that $B^T\mathcal{V}_{\mathrm{upd}}^r(v^{i_j}) $ is full rank. 
 %
 %
  % Since $V_C=S_C$, from Lemma~\ref{lemlast}, we conclude that for every $v_i$, $1\leq i\leq h$, $B^T\mathcal{V}_{\mathrm{upd}}^r(v_i)$ is full rank. 
  Accordingly,  $B^T\mathcal{V}$ is full rank, and if for some $\alpha\in \mathbb{R}^n$, $B^T\mathcal{V}\alpha=0$, we should have $\alpha=0$. In other words, $\nu=0$, contradicting the assumption. Then, the network is controllable; note that $|V_C|=n-p$. Thus, the minimum number of control nodes rendering the network controllable is $n-p$, completing the proof.
\carre

\begin{figure}[hbt]
\includegraphics[width=.24\textwidth]{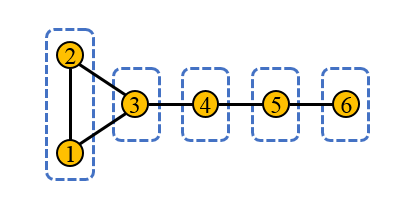}
\centering
\caption{Equitable partition in a completely uncontrollable graph.}
\label{extra}
\end{figure}

For the necessity part of Theorem \ref{thm2}, we proved that the sibling partition of a cograph is an equitable partition, where in order to break the structural symmetries, one needs to control at least all nodes of any cell except one. However, the sufficiency part of the theorem is more challenging and intricate; it holds only due to the specific structure of a cograph. To better understand the intricacy of this problem, consider a graph with a nontrivial automorphism  in Fig. \ref{extra}. One may expect that by choosing any of nodes 1 or 2 as the single control node, the controllability of this graph is ensured; while in \cite{aguilar2015graph}, it has been shown that this graph is \emph{completely uncontrollable}. Thus, breaking all symmetries in a general graph does not necessarily ensure its controllability.

 \color{black}Now, let $n_C$ be the number of different sets of control nodes with the minimum size from which a network defined on a cograph is controllable. Then, from Theorem \ref{thm2},   $n_C=\Pi_{i=1}^p m_i$. 
 
 \emph{Example:} Consider a network with dynamics (\ref{e1}) defined on the cograph in Fig. \ref{control}. By choosing a single node from cell $\mathbf{C}_1$, and two nodes from cell $\mathbf{C}_5$, the controllability of the network is ensured. For example, the network is controllable from $V_{C_1}=\{1,6,7\}$, $V_{C_2}=\{2,6,7\}$, $V_{C_3}=\{1,6,8\}$, $V_{C_4}=\{2,6,8\}$, $V_{C_5}=\{1,7,8\}$, and $V_{C_6}=\{2,7,8\}$.
 \color{black}

\subsection{Controllability of Subclasses of Cographs}

In this part, using Theorem \ref{thm2}, we derive controllability conditions for some known subclasses of cographs. 

%Let $K_1$ be a graph of a single node. Then a
A complete graph $K_n$  can be  represented as $K_n=K_1*\ldots*K_1$. By considering the corresponding cotree, one can see that $\pi_{\mathrm{sib}}(K_n)=\{\{1,\ldots,n\}\}$. Thus, a network with Laplacian dynamics (\ref{e1}) and the graph $K_n$ is controllable from at least $n-1$ nodes; a result which was established by other methods previously (e.g., see for example~\cite{zhang2011controllability,rahmani2009controllability}).

\begin{pro}
A Laplacian network (with dynamics (\ref{e1})) defined on a complete bipartite graph $K_{n_1,n_2}$ is controllable from at least $n_1+n_2-2$ control nodes. 
\end{pro} 

\emph{Proof:}
Let $V_1=\{v_1,\ldots,v_{n_1}\}$ and $V_2=\{u_1,\ldots,u_{n_2}\}$, and define $G_1=(V_1,\emptyset)$, $G_2=(V_2,\emptyset)$. Then we have  $K_{n_1,n_2}=G_1*G_2$,  implying that $\pi_{\mathrm{sib}}(K_{n_1,n_2})=\{\{v_1,\ldots,v_{n_1}\},\{u_1,\ldots,u_{n_2}\} \}$. By Theorem \ref{thm2}, the result is now immediate .\carre

In what follows, we consider an important subclass of cographs, namely threshold graphs, and as a byproduct of Theorem \ref{thm2}, we re-establish the same result of \cite{mousavi2018controllability,hsu2019minimal} by using a different approach. %extend the existing controllability results for threshold graphs~\cite{aguilar2015laplacian,hsu2016controllability}.

%\subsection{Controllability of Threshold Graphs}

%Consider the construction sequence $T^G$ associated with a threshold graph $G=(V,E)$ of size $n$. As mentioned previously, we start with a single node indexed as $1$. Then for $1<i\leq n$, if $T^G(i)=0$ (respectively, $T^G(i)=1$), in the $i$th step, a single node which is indexed with $i$ is added to the old graph through the union (respectively, join) operation. 
%As the next results show, given $T^G$, one can find the set of neighbors of any node $i\in V$.

%\begin{pro}
%In a threshold graph $G$ with the construction sequence $T^G$, for $i\in V$, if $T^G(i)=1$, then $N(i)=\{j\in V: j<i\}\cup\{j\in V: j>i,\: T^G(j)=1\}$; otherwise, $N(i)=\{j\in V: j>i,\: T^G(j)=1\}$.
%\label{pro2}
%\end{pro}
 
%\emph{Proof:} First let $T^G(i)=1$, $i\in V$. Then, the node $i$ is added to the old graph with the set of nodes $\{1,\ldots,i-1\}$ through the join operation. In other words, node $i$ is connected to all nodes $j$ for which $ j<i$. Moreover, for a node $k$ such that $k>j$, if  $T^G(k)=1$, $\{i,k\}\in E$, and if  $T^G(k)=0$, $\{i,k\}\notin E(G)$. Thus, $N(i)=\{j\in V: j<i\}\cup\{j\in V: j>i,\: T^G(j)=1\}$. On the other hand, if $T^G(i)=0$, the node $i$ is connected only to the nodes added to the graph through a join operation in step $j$ with $j>i$. In other words, $N(i)=\{j\in V: j>i,\: T^G(j)=1\}$.
%\carre 

Consider  a threshold graph $G=(V,E)$ of size $n$. Let $\mathcal{T}=(V^{\mathcal{T}}, E^{\mathcal{T}})$ be its  cotree with $v^1, \ldots, v^h$ as the internal nodes. %, where $v^h$ is the root $r$, and  $v^1$ is the furthest node from $r$. 
From the aforementioned results, one can infer the following  properties  for a threshold graph, which may not hold for a general cograph:

\begin{enumerate}
\item Any internal node of $\mathcal{T}$ has at least one child that is a leaf.
\item For any two internal nodes $v^i, v^j\in V^{\mathcal{T}}$, $\lambda_{\mathrm{upd}}^r(v^i)\neq \lambda_{\mathrm{upd}}^r(v^j)$. 
\item If $v^i$, $1\leq i\leq h$, is labeled  as 1 (respectively,  0), then for any $k\in V$ whose corresponding leaf is a child of $v^i$, we have $d(k)=\lambda_{\mathrm{upd}}^r(v^i)-1$ (respectively, $d(k)=\lambda_{\mathrm{upd}}^r(v^i)$). 
%\item For any pair of internal nodes, one is an ancestor of the other. 
\item Consider some internal node $w$ whose children--that are leaves--are indexed from $i_1$ to $i_2$. Then, the new modal matrix at $w $ is $$\scriptsize \mathcal{V}_{\mathrm{new}}(w)=
 \begin{bmatrix} 
1 &  1 & \ldots & 1\\
\vdots  &\vdots &  \ddots & \vdots\\
%1  & 1 & \ldots &  1\\
-n_1  & 1 & \ldots &1\\
0  & -n_1+1  & \ldots & 1\\
\vdots  &\vdots &  \ddots & \vdots\\
0 & 0 & \ldots & -n_2 
\end{bmatrix},
$$
where $n_1=1$ if $w$ is $v^1$, and otherwise $n_1=i_1-1$. Also, $n_2=i_2-1$. 
\end{enumerate} 

Property 1 follows from the definition of a threshold graph. Furthermore, Property 2 is inferred from Lemma \ref{lemm1}, and Properties 3 and 4  result from Theorem \ref{thmm4}. Thus, there is a direct relationship between the degree sequence of a threshold graph and its set of eigenvalues. However, a general cograph might have some internal nodes, none of whose children is a leaf.  As such, the degree sequence of a cograph does not reflect its spectrum in a transparent manner. Moreover, a general cograph might have internal nodes that are associated with the same eigenvalue, making its controllability analysis complicated. %nontrivial. 

\color{black}

 The next result  immediately follows from Theorem 1.2.4 of \cite{mahadev1995threshold}; \color{black} it states that in a threshold graph, two nodes are siblings if and only if they are of the same degree.

%Consider  a threshold graph $G=(V,E)$ of size $n$. The next result \color{Green} which immediately follow from Theorem 1.2.4 of \cite{mahadev1995threshold} \color{black} shows that in a threshold graph, two nodes are siblings if and only if they are of the same degree.
 
\begin{thm}
Given a threshold graph $G=(V,E)$, two nodes $i,j\in V$ are siblings if and only if $d(i)=d(j)$.
\label{thm3}
\end{thm}

 Now, in a threshold graph $G=(V,E)$, partition $V$ into the cells $\mathbf{C}_1,\ldots,\mathbf{C}_h$, where for any $i,j\in V$, we have $d(i)=d(j)$ if and only if for some  $1\leq k\leq h$, $i,j\in \mathbf{C}_k$.  %Let us index the cells in a way that $\mathbf{C}_i$ corresponds to a cell of leaves in $\mathcal{T}$ whose parent is the internal node $v^i$. 
  The partition $\pi_{\mathrm{deg}}=\{\mathbf{C}_1,\ldots,\mathbf{C}_h\} $ is called a \emph{degree partition}. Let $m_i=|\mathbf{C}_i|$, for $1\leq i\leq h$. %, and define $M=\max_{1\leq i\leq h} m_i$. 
   The next result that follows immediately from Theorems \ref{thm2} and  \ref{thm3} is the  same as the result of \cite{mousavi2018controllability,hsu2019minimal}  for a general threshold graph, obtained by adopting a different approach. %is an extension of the results of \cite{aguilar2015laplacian,hsu2016controllability}.
 
 \begin{cor}
For the controllability of a network defined on a threshold graph with dynamics (\ref{e1}) and the input matrix (\ref{B}), at least $n-h$ control nodes are needed, which should be chosen by selecting $m_i-1$ nodes from any cell $\mathbf{C}_i$, $1\leq i\leq h$, of its degree partition.
% Let $\pi_{\mathrm{deg}}=\{\mathbf{C}_1,\ldots,\mathbf{C}_h\}$ be a degree partition in a connected threshold graph $G$, where $m_i=|\mathbf{C}_i|$, for $1\leq i\leq h$. Consider a network defined on $G$ with dynamics (\ref{e1}) and the input matrix (\ref{B}).  Then, at least $n-p$ control nodes are needed to render the network controllable, which should be chosen by selecting $m_i-1$ nodes from any cell $\mathbf{C}_i$.   
 \end{cor}

\section{Conclusions}

%Strong structural controllability of the zero and nonzero modes of networks have been examined in this paper.
%, which can be respectively interpreted as the controllability of the null spaces of matrices carrying the structure of the network and their orthogonal complements. 

%In this work, a new notion of strong zero forcing set is introduced which can be used in finding an upper bound on the maximum geometric multiplicity of nonzero eigenvalues for a family of pattern matrices. 
%Moreover,
In this work, we characterized the controllability of Laplacian networks defined over cographs in terms of certain graph-theoretic conditions. 
These characterizations are built upon the intricate correspondance between the inherent structural modularity of cographs, with respect to join and union operation, and its modal properties.
Moreover, we used the proposed framework 
to provide a procedure for selecting a set of control nodes guaranteeing the controllability of cograph networks. In particular, we demonstrated that the minimum number of control nodes rendering a cograph controllable is the difference between its size and the number of cells of its sibling partition. It was also revealed that the larger a cell of sibling nodes, the larger the multiplicity of one of the eigenvalues associated with the Laplacian matrix; such multiplicities are often associated with higher degrees of symmetry in the network. 
%Hence, at least intuitively, more nodes are required as control nodes to break such symmetries for network controllability.
% 
We then applied our results to certain subclasses of cographs such as threshold graphs, and by adopting a different approach,  presented conditions that ensure their controllability. As one of the future research topics, in order to minimize the number of independent controllers of a cograph, one can consider the more general case where the input matrix $B$ is binary with more than one nonzero entry at each column. Furthermore, controllability analysis of weighted or directed cographs can be taken into account. %, extending previous results reported in the literature.
\balance
{
\section*{Acknowledgements}
The authors thank Dr. Cesar Aguilar for insightful discussions  on network controllability and pointing out a streamlined approach to  derive \color{black} the necessary controllability condition for cographs using their inherent symmetry.}
%
%{\color{blue} for references it might be better if you consistently use the first initials only, e.g., M. Mesbahi; if you want to use full first names, then make them consistent, e.g., check [30], not sure about the format for [29].}
%
%\bibliographystyle{plain}        % Include this if you use bibtex 
%\bibliographystyle{elsarticle-harv}
%\bibliographystyle{elsarticle-num-names}

\bibliographystyle{IEEEtran}
\bibliography{library}

\end{document}